\documentclass{article}
\usepackage[T2A]{fontenc}
\usepackage[utf8]{inputenc}
\usepackage[english,russian]{babel}
\usepackage[tbtags]{amsmath}
\usepackage{amsfonts,amssymb,mathrsfs,amscd}
\usepackage[hyper]{msb-a}

\makeatletter
\gdef\No{{\select@language{russian}\textnumero}}
\makeatother

\JournalName{}
\numberwithin{equation}{section}

\def \N {{\mathbf {N}}}

\def \T {{\mathbf  T}}

\begin{document}


\title{An ergodic automorphism $\bf T$ with  singular spectrum  of $\bf T^{\otimes n}$ and  Lebesgue  one of  $\bf T^{\otimes (n+1)}$  }

\author[В.\,В.~Рыжиков]{V.\,V.~Ryzhikov} %
\address{Московский государственный университет}
\email{vryzh@mail.ru}

\date{2024}
\udk{517.987}

\maketitle

\begin{fulltext}
\begin{abstract}{

For any natural $n$, and real  $\alpha\geq 0$ we construct 
a Sidon automorphism $T$ such that its tensor powers $T^{\otimes n}$
 have singular spectrum if $n\leq 1+\alpha /2$, and  Lebesgue spectrum if $n\, > 1+\alpha/2$. Moreover, the power $T^{\otimes n}$ is conservative if $n\leq 1+\alpha$ and dissipative if
$n>1+\alpha$.

}
\end{abstract}
\markright{ }



\vspace{2mm}
\bf Theorem 1. \it For any $n$ 
there is an ergodic automorphism $T$ of a space with sigma-finite measure whose power $T^{\otimes 2n}$ is conservative,  spectrum of $T^{\otimes (n+1)}$ is Lebesgue, and spectrum
of  $T^{\otimes n}$ is singular.
\rm

\vspace{2mm}
Thus, ergodic transformations are sources of measures $\sigma$ with a given singularity order $n$:   the convolution powers  $\sigma^{\ast m}$ are singular as $m\leq n$ and are  absolutely continuous for  $m > n$.

\vspace{2mm}
Specially constructed transformations of rank 1 were  used, for example, to  answer \it spectral \rm questions  of Kolmogorov, Rokhlin, Gordin, Bergelson, Oseledets (see \cite{23}).
Examples of automorphisms satisfying the spectral properties specified in Theorem 1 are
 of rank 1. Let us recall their definition.

\vspace{2mm}
\bf Constructions of rank one. \rm
We fix a natural number $h_1$, a sequence $r_j$ (the number of columns into which the tower of stage $j$ is cut) and a sequence of integer vectors (spacer parameters)
$$ \bar s_j=(s_j(1), s_j(2),\dots, s_j(r_j-1),s_j(r_j)).$$

Let a system of disjoint half-intervals be defined at step $j$
$$E_j, TE_j, T^2E_j,\dots, T^{h_j-1}E_j,$$
and on half-intervals $E_j, TE_j, \dots, T^{h_j-2}E_j$
the transformation $T$ is a parallel transfer. Such a set of half-intervals is called a tower of stage $j$; their union is denoted by $X_j$ and is also called a tower.

Let us represent $E_j$ as a disjoint union of $r_j$ half-intervals
$$E_j^1,E_j^2E_j^3,\dots E_j^{r_j}$$ of the same length.
For each $i=1,2,\dots, r_j$ we define the column $X_{i,j}$ as the union of intervals
$$E_j^i, TE_j^i ,T^2 E_j^i,\dots, T^{h_j-1}E_j^i.$$

To each column $X_{i,j}$ we add $s_j(i)$ of disjoint half-intervals of the same measure as $E_j^i$, obtaining a set
$$E_j^i, TE_j^i, T^2 E_j^i,\dots, T^{h_j-1}E_j^i, T^{h_j}E_j^i, T^{h_j+1}E_j^i , \dots, T^{h_j+s_j(i)-1}E_j^i$$
(all these sets do not intersect).
Denoting $E_{j+1}= E^1_j$, for $i<r_j$ we set
$$T^{h_j+s_j(i)}E_j^i = E_j^{i+1}.$$
 The set of superstructured columns is from now on considered as a tower of stage $j+1$, consisting of half-intervals
$$E_{j+1}, TE_{j+1}, T^2 E_{j+1},\dots, T^{h_{j+1}-1}E_{j+1},$$
where
 $$ h_{j+1}+1 =(h_j+1)r_j +\sum_{i=1}^{r_j}s_j(i).$$

As a result we get an invertible transformation $T:X\to X$  defined on the space $X=\cup_j X_j$, it preserves the standard Lebesgue measure on $X$.

An automorphism of rank one is ergodic and has a simple  spectrum. It is known that indicators of intervals appearing in the description of  constructions $T$ are cyclic vectors for our constructions  $T$.

\bf Sidon automorphisms of class $\bf C(\alpha)$. \rm Rank one construction $T$ is called Sidon if it has the following property: \it $T^mX_j$ may overlap   only with  one of the columns $X_{i,j}$  as $h_{j}<m\leq h_{j+1}$. 
 \rm 
\vspace{2mm}

 We say that a Sidon construction $T$ 
belongs \it to class $\bf C(\alpha)$, \rm if for some sequence
$j(k)\to\infty$ for any $\delta>0$
$$\sum_k \frac 1 {r_{j(k)}^\delta}<\infty$$
$$r_j=r_{j(k)}, \ j(k)\leq j<j(k+1), $$
$$j(k+1) -j(k)=[r_{j(k)}^\alpha].$$

\vspace{2mm}
\bf Theorem 2. \it
Let the Sidon automorphism $T$ belong to the class
$\bf C(\alpha)$, $\alpha\geq 0$.
If $n\leq 1+\alpha /2$, the power $T^{\otimes n}$ has
singular spectrum. Its spectrum is  absolutely continuous as $n> 1+\alpha/2$.
If $n\leq 1+\alpha$, then  $T^{\otimes n}$ is conservative, 
and  dissipative as $n> 1+\alpha$.
\rm

\vspace{2mm}
\bf The  construction $\bf T$ whose power $\bf T^{\otimes 20}$ is conservative, 

spectrum of $\bf T^{\otimes 11}$ is Lebesgue, 

and spectrum of $\bf T^{\otimes 10}$ is  singular.\rm

\vspace{3mm}
Let's put
$j(k)= \sum_{m=1}^k (m!)^{20}$  
and  define the  parameters of $T$:

$r_{j}:= (k+1)!, \j(k)\leq j <j(k+1)$,

$s_j(i):=10^ih_j, \ 1 \leq i \leq r_j,$

where $h_1=10$, $h_{j+1}= r_jh_j +\sum_{i=1}{r_j} s_j(i)$.


\newpage
We  formulate below statements from which the stated results are derived.

\vspace{2mm}
\bf Theorem 2.1. \it Let a Sidon construction $T$ satisfy
$$\sum_{i=1}^{\infty }\left(\frac 1 {r_j}\right)^{d-1} =\infty,$$
then  $T^{\otimes d}$ is  conservative: there is no wandering set 
 of positive measure  in $X^d$.
\rm

The above assertion is proved in more general case in \cite{LS}, theorem 3.2.

\vspace{3mm}
\bf Theorem 2.2. \it
If for a Sidon construction $T$ we have 
$$\sum_{i=1}^{\infty }\left(\frac 1 {r_j}\right)^{2d-2}<\infty,$$
then spectrum of $T^{\otimes d}$ is absolutely continuous (and it is guaranteed to be Lebesgue under some additional conditions on $r_j$).

\rm

\vspace{3mm}
\bf Theorem 2.3. \it For a Sidon construction $T$ of class $\bf C(\alpha)$ the condition 
 $$\sum_{k=1}^{\infty }\left(\frac 1 {r_{j(k)}}\right)^{2d-2-\alpha}=\infty$$ 
implies the singularity of  spectrum of   $T^{\otimes d}$.\rm

\vspace{2mm}
Proofs of these statements use modified ideas from \cite{24}.

\vspace{2mm}

 It would be interesting to find out the ergodic properties of conservative powers with an absolutely continuous spectrum. It is also of  interest to clarify the spectral properties of ergodic tensor powers of transformations from the work \cite{LS}.

\section{Эргодический  автоморфизм $\bf T$, для которого  спектр у $\bf T^{\otimes (n+1)}$  лебеговский, а спектр степени  $\bf T^{\otimes n}$  сингулярный} 

Для всякого действительного числа $\alpha\geq 0$ и натурального $n$ мы предъявляем 
 автоморфизм $T$ такой, что  тензорные степени $T^{\otimes n}$ обладают сингулярным спектром  при $n\leq 1+\alpha /2$ и лебеговским спектром  при $n\, > 1+\alpha/2$.  Причем степени $T^{\otimes n}$ консервативны в случае $n\leq 1+\alpha$ и диссипативны при 
$n>1+\alpha$.

Унитарный оператор $S$ на сепарабельном гильбертовом пространстве, 
как известно,  изоморфен оператору $V$, $Vf(t,n)=zf(t,n),$ 
который действует в пространстве
$L_2(\T\times \N, \sigma)$, где $\T$  -- единичная окружность в комплексной плоскости,
$\sigma$ -- нормированная борелевская мера на $\T\times \N$. 

Если мера $\sigma$ сингулярна, то говорят, что спектр оператора $S$ сингулярен. Если проекция меры $\sigma$ на  $\T$ эквивалентна мере Лебега на $\T$, то спектр такого  $S$  называется лебеговским.

\vspace{2mm}
\bf Теорема 1. \it Для всякого $n$
найдется эргодический автоморфизм $T$ пространства с сигма-конечной мерой, чья  степень  $T^{\otimes 2n}$ консервативна, спектр $T^{\otimes (n+1)}$  лебеговский, а спектр  $T^{\otimes n}$  сингулярен.
\rm

\vspace{2mm}
Теорема доказана для $n=1$ в  \cite {24}.
Специально построенные преобразования ранга 1  использовались, например,  для ответов на вопросы Колмогорова, Рохлина, Гордина, Бергельсона, Оселедца  (см. \cite{23}).
Примеры автоморфизмов, удовлетворяющих спектральным свойствам, указанным в теореме,
также удобно искать в классе преобразований ранга 1. Напомним их определение.

\vspace{2mm}
\bf   Конструкции ранга один. \rm
Фиксируем натуральное число $h_1$, последовательность  $r_j$ (число колонн, на которые  разрезается башня этапа $j$) и  последовательность целочисленных векторов (параметров надстроек)   
$$ \bar s_j=(s_j(1), s_j(2),\dots, s_j(r_j-1),s_j(r_j)).$$  

Пусть на шаге $j$  определена система   непересекающихся полуинтервалов 
$$E_j, TE_j, T^2E_j,\dots, T^{h_j-1}E_j,$$
причем на полуинтервалах $E_j, TE_j, \dots, T^{h_j-2}E_j$
пребразование $T$ является параллельным переносом. Такой набор   полуинтервалов  называется башней этапа $j$, их объединение обозначается через $X_j$ и также называется башней.

Представим   $E_j$ как дизъюнктное объединение  $r_j$ полуинтервалов 
$$E_j^1,E_j^2E_j^3,\dots E_j^{r_j}$$ одинаковой длины.  
Для  каждого $i=1,2,\dots, r_j$ определим  колонну $X_{i,j}$ как объединение интервалов  
$$E_j^i, TE_j^i ,T^2 E_j^i,\dots, T^{h_j-1}E_j^i.$$

К каждой  колонне $X_{i,j}$ добавим  $s_j(i)$  непересекающихся полуинтервалов  той же меры, что у $E_j^i$, получая набор  
$$E_j^i, TE_j^i, T^2 E_j^i,\dots, T^{h_j-1}E_j^i, T^{h_j}E_j^i, T^{h_j+1}E_j^i, \dots, T^{h_j+s_j(i)-1}E_j^i$$
(все эти множества  не пересекаются).
Обозначив $E_{j+1}= E^1_j$, для   $i<r_j$ положим 
$$T^{h_j+s_j(i)}E_j^i = E_j^{i+1}.$$
 Набор надстроеных колонн с этого момента  рассматривается как   башня  этапа $j+1$,  состоящая из полуинтервалов  
$$E_{j+1}, TE_{j+1}, T^2 E_{j+1},\dots, T^{h_{j+1}-1}E_{j+1},$$
где  
 $$ h_{j+1}+1 =(h_j+1)r_j +\sum_{i=1}^{r_j}s_j(i).$$

Частичное определение преобразования $T$ на этапе $j$ сохраняется на всех следующих этапах. В итоге на пространстве  $X=\cup_j X_j$ определено  обратимое преобразование $T:X\to X$, сохраняющее  стандартную меру Лебега на $X$.

Преобразование (автоморфизм) и индуцированный им унитарный оператор в статье обозначаются одинаково.
Автоморфизм ранга один эргодичен,   имеет простой (однократный) спектр.  Известно, что линейные комбинации индикаторов  интервалов, фигурирующих  в описании конструкций  $T$, являются циклическими векторами для оператора $T$.   

 \bf   Сидоновские автоморфизмы класса $C(\alpha)$. \rm Пусть конструкция $T$ ранга один обладает следующим свойством: \it пересечение 
$T^mX_j$ при   $h_{j}<m\leq h_{j+1}$  может содержаться 
только в одной  из колонн  $X_{i,j}$ башни $X_j$. \rm Такая конструкция 
называется  \it сидоновской. \rm
\vspace{2mm}

  Говорим, что  сидоновская конструкция  $T$ ранга один
принадлежит \it классу $C(\alpha)$, \rm если для некоторой последовательности 
$j(k)\to\infty$ для любого $\delta>0$ имеем
$$\sum_k  \frac 1 {r_{j(k)}^\delta}<\infty$$  
$$r_j=r_{j(k)}, \  j(k)\leq j<j(k+1), $$
$$j(k+1) -j(k)=[r_{j(k)}^\alpha].$$

\vspace{2mm}
\bf Теорема 2. \it 
Пусть сидоновский автоморфизм $T$ принадлежит классу
$C(\alpha)$, $\alpha\geq 0$. 
Если $n\leq 1+\alpha /2$,   степень $T^{\otimes n}$ обладает 
сингулярным спектром,  при $n> 1+\alpha/2$    абсолютно непрерывным. 
Если  $n\leq 1+\alpha$, то   степень $T^{\otimes n}$ является консервативной, а при $n> 1+\alpha$ диссипативной.
\rm

\vspace{2mm}
\bf Пример конструкции  $\bf T$, для которой выполнено:  

степень  $\bf T^{\otimes 20}$ консервативна, 

спектр $\bf T^{\otimes 11}$  лебеговский, 

спектр  $\bf T^{\otimes 10}$  сингулярен.\rm

Положив 
$j(k)= \sum_{m=1}^k (m!)^{20}$,

определим параметры конструкции:

$r_{j}:= (k+1)!, \ j(k)\leq j <j(k+1)$,

$s_j(i):=10^ih_j,  \ 1 \leq  i \leq r_j,$  

где  $h_1=10$,  $h_{j+1}= r_jh_j +\sum_{i=1}{r_j} s_j(i)$.


Некоторые вспомогательные утверждения.

\vspace{2mm}
\bf Теорема 2.1.   \it Пусть для сидоновской конструкции $T$ выполнено
$$\sum_{i=1}^\infty\left(\frac 1 {r_j}\right)^{d-1} =\infty,$$ 
тогда степень $T^{\otimes d}$ является консервативным преобразованием. 
\rm

\vspace{2mm}
Вытекает непосредственно из  \cite{LS}, теорема 3.2.  

\vspace{2mm}
\bf Теорема 2.2.   \it 
Если для сидоновской конструкции $T$  сходится ряд
$$\sum_{i=1}^{\infty }\left(\frac 1 {r_j}\right)^{2d-2},$$ 
то спектральная мера степени $T^{\otimes d}$ абсолютно непрерывна.

\rm
 
\vspace{2mm}
\bf Теорема 2.3.   \it Для сидоновской конструкции $T$ класс $C(\alpha)$ cпектр
степени   $T^{\otimes d}$ сингулярен, если выполнено условие 
 $$\sum_{k=1}^{\infty }\left(\frac 1 {r_{j(k)}}\right)^{2d-2-\alpha}=\infty.$$ \rm  

\vspace{2mm}
Доказательства утверждений используют модификацию методов  работы \cite{24}.

\vspace{2mm}
 Интересно  выяснить эргодические (метрические) свойства консервативных степенией с абсолютно непрерывным спектром.  Также представляет  интерес выяснение спектральных свойств эргодических тензорных степеней преобразований из  работы \cite{LS}.

\end{fulltext}


\begin{thebibliography}{99}

\bibitem{23}В.В. Рыжиков, Спектры самоподобных эргодических действий, Матем. заметки, 113:2 (2023),  273-282; 

V.V. Ryzhikov, Spectra of Self-Similar Ergodic Actions, Math. Notes, 113:2 (2023), 274-281



\bibitem{24}В.В. Рыжиков, 
Полиномиальная жесткость и спектр сидоновских автоморфизмов. Матем. сб., 215:7 (2024)7;
  
V.V. Ryzhikov, Polynomial rigidity and spectrum of Sidon automorphisms. Sb. Math., 215:7 (2024)

\bibitem{LS}I. Loh,  C.E. Silva, Strict doubly ergodic infinite transformations, 
Dyn. Syst. 32, No. 4, 519-543 (2017). 


\end{thebibliography}
\end{document}